\documentclass[12pt]{amsart}
\newcommand{\bC}{{\mathbb C}}

\newcommand{\bP}{{\mathbb P}}
\newcommand{\bQ}{{\mathbb Q}}

\newcommand{\cC}{{\mathcal C}}
\newcommand{\cE}{{\mathcal E}}
\newcommand{\cF}{{\mathcal F}}

\newcommand{\cH}{{\mathcal H}}
\newcommand{\cI}{{\mathcal I}}
\newcommand{\cJ}{{\mathcal J}}
\newcommand{\cK}{{\mathcal K}}
\newcommand{\cL}{{\mathcal L}}
\newcommand{\cO}{{\mathcal O}}

\newcommand{\cQ}{{\mathcal Q}}

\newcommand{\oX}{\overline{X}}

\newcommand{\tC}{\widetilde{C}}

\newcommand{\tkappa}{\widetilde{\kappa}}

\newcommand{\tV}{\widetilde{V}}
\newcommand{\tX}{\widetilde{X}}
\newcommand{\tZ}{\widetilde{Z}}
\newcommand{\ra}{\rightarrow}
\newcommand{\lra}{\longrightarrow}

\newcommand{\inj}{\hookrightarrow}
\newcommand{\surj}{-\hspace{-4pt}\rightarrow\hspace{-19pt}{\rightarrow}\hspace{3pt}}

\newcommand{\T}{\Theta}

\newcommand{\Sg}{\mathfrak{S}}

\newtheorem{theorem}{Theorem}[section]
\newtheorem{lemma}[theorem]{Lemma}

\numberwithin{equation}{section}

\begin{document}

\title{Deforming curves in jacobians to non-jacobians II: \\
curves in $C^{ (e) }, 3\leq e\leq g-3$}

\author{E. Izadi}

\address{Department of Mathematics, Boyd
Graduate Studies Research Center, University of Georgia, Athens, GA
30602-7403, USA}

\email{izadi@math.uga.edu}

\thanks{This material is based upon work partially supported by the
National Science Foundation under Grant No. DMS-0071795. Any opinions,
findings and conclusions or recomendations expressed in this material
are those of the author and do not necessarily reflect the views of
the National Science Foundation}



\maketitle

\section*{Introduction}

This is a second paper where we introduce deformation theory
methods which can be applied to finding curves in families of
principally polarized abelian varieties ({\em ppav}) containing
jacobians. One of our motivations for finding interesting and
computationally tractable curves in ppav is to solve the Hodge
conjecture for the primitive cohomology of the theta divisor which we
explain below. For other motivations we refer to the prequel
\cite{I13} to this paper.

Let $(A ,\T)$ be a ppav over $\bC$ of dimension $g\geq 4$ such that
$\T$ is smooth. Since any abelian variety (over $\bC$) is isogenous to
such an abelian variety, the Hodge conjectures for arbitrary abelian
varieties are equivalent to the Hodge conjectures for principally
polarized abelian varieties with smooth theta divisors.

The primitive part $K(\T, \bQ)$ of the cohomology of $\T$ can be
defined as the kernel of the map $H^{g-1}(\T, \bQ) \lra H^{g+1}(A,
\bQ)$ obtained by Poincar\'e Duality from push-forward on
homology. The space $K(\T ,\bQ)$ defines a Hodge substructure of the
cohomology of $\T$ of level $g-3$ (see page 562 of \cite{I4}; the
proof there works also for $g > 4$). The generalized Hodge conjecture
then predicts that there is a family of curves in $\T$ such that
$K(\T, \bQ)$ is contained in the image of its Abel-Jacobi map. The
Abel-Jacobi map for a family of curves can be defined as follows.

Let $\cC \ra S$ be a family of curves with $S$ smooth,
complete and irreducible of dimension $d$ such that there is a diagram
\[
\begin{array}{rcc}
\cC & \stackrel{q}{\lra} & \T \\
p \downarrow & & \\
S & & .
\end{array}
\]
The Abel-Jacobi map for this family of curves is by definition $AJ :=
q_*p^* : H^{2d-(g-3)}(S, \bQ) \ra H^{g-1}(\T, \bQ)$. The image of
the Abel-Jacobi map defines a Hodge substructure of level $\leq g-3$
of the cohomology of $\T$.

For abelian fourfolds one interesting family of curves is the family
of Prym-embedded curves in $\T$ and it is proved in \cite{I4} that it
does give a solution to the Hodge conjecture for $K (\T, \bQ)$. In
dimension $\geq 6$ there are no known families of interesting
curves in the theta divisor of a general ppav.

Let us briefly explain our methods, similar to \cite{I13}. After
identifying $JC = Pic^0 C$ with $A := Pic^{ g-1 } C$ by tensoring with a
fixed invertible sheaf of degree $g-1$, Riemann's theta divisor is
\[
\Theta :=\{\cL\in Pic^{ g-1 } C : h^0 (\cL ) > 0\}.
\]
Consider a subvariety $X$ of $A$ contained in ``many'' translates
$\T_a$ of $\T$. As in \cite{I13}, for each such translate $\T_a$, we have a map
\[
\nu_a : H^1 (T_A )\lra H^1 (\cI_{ Z_{ g-1 }} (\Theta_a ) |_X),
\]
obtained from Green's exact sequence (\cite{green84}, see Section
\ref{sectrans} below) which factors through the first order obstruction map
\[
\nu : H^1 (T_A )\lra H^1 (N_{ X/JC })
\]
where $N_{ X/A }$ is the normal sheaf to $X$ in $A$(see Section
\ref{sectobs}). Hence, if $\nu_a (\eta)$ is not zero for some $a$, so
is $\nu(\eta)$.

The main difference between the method used here and that of
\cite{I13} is in Section \ref{sectker1} below which is more
difficult for $e >2$ and is where we need some assumptions of
genericity on $X$.

We apply the above to families of curves in jacobians which are
natural generalizations of Prym-embedded curves in tetragonal
jacobians. More precisely, let $C$ be a non-hyperelliptic curve of genus
$g$ with a $g^1_d$ (a pencil of degree $d$). Define
\[
X_e (g^1_d) :=\{ D_e : \exists D\in C^{ (d-e) } \hbox{ such that } D_e
+ D\in g^1_d\}\subset C^{ (e) }
\]
where $2\leq e\leq d$ and $C^{ (e) }$ is the $e$-th symmetric power of
$C$ (see Section \ref{defX} below for the precise definition).  We map
$X := X_e (g^1_d)$ and $C^{ (e) }$ to $C^{ (g-1) }$ and then to $A$ by
adding a fixed divisor $q :=\sum_{ i=1 }^{ g-1-e } q_i$. If $d\geq
e+1$, the map is non-constant on $X$. We call $W_e+q$ the image of
$C^{ (e) }$ in $A$ via this map. Given a one-parameter infinitesimal
deformation of the jacobian of $C$ normal to the jacobian locus
$\cJ_g$ we ask when the curve $X$ deforms with it. Let $Z_{ g-1
}\subset C^{ (g-1) }$ be the locus where the map
\[\begin{array}{cccl}
\rho : & C^{ (g-1) } & \lra & \Theta\subset A \\
 & D & \longmapsto & \cO_C (D)
\end{array}
\]
fails to be an isomorphism and let $Z_q\subset
C^{ (e) }$ be the locus defined by the exactness of the sequence
\[
N_{ W_e+q/A } |_{ C^{ (e) }}\lra N_{\Theta /A} |_{ C^{ (e) }} =\cO_{
C^{ (e) }} (\Theta )\lra\cO_{ Z_q } (\T )\lra 0.
\]
We prove the following
\begin{theorem}\label{mainthm}
  Assume $3\leq e\leq g-3$, for all $e'\leq e$ the curve $X_{e'}
(g^1_d)$ is irreducible and reduced and the set of $q$ for which
$X\cap Z_{ g-1 }\neq X\cap Z_q$ has dimension at most $g-e-3$. If $X_e
(g^1_d)$ deforms out of $\cJ_g$ then
\begin{itemize}
\item either $h^0 (g^1_d) =e$ and $d=2e$
\item or $h^0 (g^1_d ) =e+1$ and $d=2e+1$.
\end{itemize}
\end{theorem}
By Appendix \ref{generic} below, for $(C, L)$ in a non-empty open subset of
the (irreducible) Hurwitz scheme of smooth curves with maps of degree
$d$ to $\bP^1$ (with simple ramification) the hypotheses of the
theorem are satisfied. In case $e=2$, we proved this result in
\cite{I13} without the assumptions of genericity. We expect that for
$e > 2$ the result will still hold for reducible curves but
non-reduced curves might deform in directions which are contained
in the intersection of the spans, in $S^2 H^1 (\cO_A )\subset H^1
(\cO_A )^{\otimes 2}\cong H^1 (T_A )$ of the divisors parametrized by the
curve $X$ (this intersection is empty for reduced curves but could be
non-empty for non-reduced curves).

In the case $e=2$, $d=4$, the curve $X$ is a Prym-embedded curve (see
Recillas \cite{recillas74}), hence deforms out of $\cJ_g$ into the
locus of Prym varieties.


As explained in Section \ref{sectconsq}, Theorem \ref{mainthm} shows
that when $3\leq e\leq g-3$ the most interesting cases in which $X$
will possibly deform out of $\cJ_g$ are those in which $C$ is
bielliptic, $e$ any integer between $3$ and $g-3$ or $C$ is any curve
of genus $g$ between $6$ and $10$, $e=3$ and $L\subset g^2_6$. In both
these cases, it is likely that the curve $X$ will deform out of
$\cJ_g$. Although jacobians of bielliptic curves form a subvariety of
dimension $2g-2$ of $\cJ_g$, the curves obtained from bielliptic
curves will likely deform to large families of ppav: such a situation
is analogous to the case of tetragonal jacobians of dimension $\geq 7$
where the curve $X_2 (g^1_4 )$ deforms to a general Prym but does NOT
deform to a general jacobian.

So we have some families of curves (including any $X$ with $e=g-2$) which
could possibly deform to non-jacobians. We need a different approach
to prove that higher order obstructions to deformations vanish: this
will be presented in detail in the forth-coming paper \cite{I16} and
the idea behind it is the following. For each $\T_a$ containing $X$,
one has the map of cohomology groups of normal sheaves
\[
H^1 (N_{ X/JC })\lra H^1 (N_{\T_a /JC } |_X ) = H^1 (\cO_X (\T_a ))
\]
whose kernel contains all the obstructions to the deformations of $X$
since we will only consider algebraizable deformations of $JC$ for
which the obstructions to deforming $\T_a$ vanish. If one can prove
that the intersection of these kernels is the image of the first order
algebraizable deformations of $JC$, i.e., the image of $S^2 H^1 (\cO_C
)\subset H^1 (T_{ JC})$, it will follow that the only obstructions to
deforming $X$ with $JC$ are the first order obstructions.

Finally, we would like to mention that from curves one can obtain
higher-dimensional subvarieties of an abelian variety. For a
discussion of this we refer the reader to \cite{I14}.

Plan of the paper: In Section \ref{defX} we define the curves $X_e
(L)$ via their ideals for which we write down a concrete workable
resolution. We compute their genus, define their maps to $A$ and
prove a useful lemma about divisors parametrized by $X_{ e+1 }(L)$ and
$X_{ e+2 }(L)$. In Section \ref{sectobs} we define the first order
obstruction map $\nu_e : S^2 H^1 (\cO_C )\ra H^1 (N_{ W_e/A} |_X)$
which we use to prove Theorem \ref{mainthm}. In Section \ref{sectrans}
we compute the translates $\T_a$ of $\T$ containing $X$ and show how
we can ``replace'' $\nu_e$ by the collection of maps $S^2 H^1
(\cO_C)\ra H^1 (\cI_{Z_{ g-1 }\cap X_{ -a }} (\T))$. Our method will
consist in finding when these maps can have non-trivial kernels. In
Section \ref{sectkertot} we decompose these maps into compositions of
$S^2 H^1 (\cO_C )\ra H^0 (\cO_{ Z_{ g-1 }\cap X_{-a}} (\Theta ))$ and
$H^0 (\cO_{ Z_{ g-1 }\cap X_{-a}} (\Theta ))\ra H^1 (\cI_{ Z_{ g-1
}\cap X_{ -a }} (\Theta ))$ which we then analyze separately in
Sections \ref{sectker1} and \ref{sectker2} respectively. In Section
\ref{sectker1} we prove that for any $\eta\in S^2 H^1 (\cO_C)\setminus
H^1 (T_C)$, there exists a translate $\T_a$ of $\T$ containing $X$
such that $\eta$ is {\em not} in the kernel of $S^2 H^1 (\cO_C )\ra
H^0 (\cO_{ Z_{ g-1 }\cap X_{-a}} (\Theta ))$. In Section
\ref{sectker2} we prove that for ``almost all'' $\T_a$ containing $X$, the
coboundary map $H^0 (\cO_{ Z_{ g-1 }\cap X_{-a}} (\Theta ))\ra H^1
(\cI_{ Z_{ g-1 }\cap X_{ -a }} (\Theta ))$ is injective unless $d=2e$
and $h^0 (L) =e$ or $d=2e+1$ and $h^0(L) = e+1$ which proves Theorem
\ref{mainthm}. In Section \ref{sectconsq} we describe the consequences
of Theorem \ref{mainthm}. Finally, we gather some useful technical
results in the Appendix.

\section*{Notation and Conventions}

We will denote linear equivalence of divisors by $\sim$.

For any divisor or coherent sheaf $D$ on a scheme $X$, denote by $h^i
(D)$ the dimension of the cohomology group $H^i (D) = H^i (X, D)$. For
any subscheme $Y$ of $X$, we will denote by $\cI_{Y/X}$ the ideal sheaf
of $Y$ in $X$ and by $N_{ Y/X }$ the normal sheaf of $Y$ in $X$. When
there is no ambiguity we drop the subscript ${}_X$ from $\cI_{Y/X}$ or
$N_{ Y/X }$. The tangent sheaf of $X$ will be denoted by $T_X :=\cH om
(\Omega^1_X ,\cO_X )$ and the dualizing sheaf of $X$ by $\omega_X$. By
the genericity of any property on $X$ we mean genericity on every
irreducible component.

We let $C$ be a smooth nonhyperelliptic curve of genus $g$ over the
field $\mathbb C$ of complex numbers. For any positive integer $n$,
denote by $C^{ (n) }$ the $n$-th symmetric power of $C$. Note that
$C^{ (n) }$ parametrizes the effective divisors of degree $n$ on $C$.

We denote by $K$ an arbitrary canonical divisor on $C$. Since $C$ is
not hyperelliptic, its canonical map $C\ra |K|^*$ is an embedding and
throughout this paper we identify $C$ with its canonical image. For a
divisor $D$ on $C$, we denote by $\langle D\rangle$ its span in $| K
|^* =\bP H^0 (\omega_C )^* =\bP H^1 (\cO_C)$.

Since we will mostly work with the Picard group $Pic^{ g-1 }C$ of
invertible sheaves of degree $g-1$ on $C$, we put $A := Pic^{g-1}
C$. Let $\Theta$ denote the natural theta divisor of $A$, i.e.,
\[
\Theta :=\{\cL\in A : h^0 (\cL) > 0\}\; .
\]
The multiplicity of $\Theta$ at $\cL\in\Theta$ is $h^0 (\cL)$
(\cite{ACGH} Chapter VI p. 226). So the singular locus of $\Theta$ is
\[
Sing (\Theta ):=\{\cL\in A : h^0 (\cL )\geq 2\}\; .
\]
There is a map
\[
\begin{array}{rcl}
Sing_2 (\Theta ) & \lra & |I_2 (C)|\\
\cL & \longmapsto & Q (\cL ) :=\cup_{ D\in |\cL | }\langle D\rangle
\end{array}
\]
where $Sing_2 (\Theta )$ is the locus of points of order $2$ on
$\Theta$ and $| I_2 (C) |$ is the linear system of quadrics containing
the canonical curve $C$. This map is equal to the map sending $\cL$ to
the (quadric) tangent cone to $\Theta$ at $\cL$ and its image $\cQ$
generates $| I_2 (C) |$ (see \cite{green84} and
\cite{smithvarley90}). Any $Q (\cL )\in\cQ$ has rank $\leq 4$. The
singular locus of $Q (\cL)$ cuts $C$ in the sum of the base divisors
of $|\cL |$ and $ |\omega_C\otimes\cL^{ -1 } |$. The rulings of $Q$
cut the divisors of the moving parts of $|\cL |$ and $
|\omega_C\otimes\cL^{ -1 } |$ on $C$ (see \cite{andreottimayer67}).

For any divisor or invertible sheaf $a$ of degree $0$ and any
subscheme $Y$ of $A$, we let $Y_a$ or $Y+a$ denote the translate of $Y$ by
$a$. By a $g^r_d$ we mean a (not necessarily complete) linear system
of degree $d$ and dimension $r$. We call $W^r_d$ the subvariety of
$Pic^d C$ parametrizing invertible sheaves $\cL$ with $h^0 (\cL) > r$.

For any effective divisor $E$ of degree $e$ on $C$ and any positive
integer $n\geq e$, let $C^{( n-e) }_E\subset C^{ (n )}$ be the image
of $C^{ (n-e )}$ in $C^{ (n) }$ by the morphism $D\mapsto D+E$. For
any divisor $E =\sum_{i=1}^{ r} n_i t_i$ on $C$, let $C^{ div }_E$
denote the divisor $\sum_{ i=1 }^{ r} n_i C_{t_i}^{ (n-1) }$ on $C^{
(n) }$. For a linear system $L$ on $C$, we denote by $C^{ div }_L$ any
divisor $C^{ div }_E\subset C^{ (n) }$ with $E\in L$.

By infinitesimal deformation we always mean {\em flat} first order
infinitesimal deformation.

\section{The curve $X := X_e (g^1_d)$ and its useful properties}
\label{defX}

Suppose $2\leq e\leq g-1$ and let $L$ be a pencil of degree $d\geq
e+2$ on $C$. We would like to define a curve $X$ whose underlying set
will be
\[
\{ D_e : \exists D\in C^{ (d-e) } \hbox{ such that } D_e + D\in
L\}.
\]
If $L$ contains reduced divisors, then $X$ is reduced and can be
defined by the above set. If $L$ does not contain reduced divisors
then we need to define a scheme structure on $X$. Although we suppose
$X$ integral in this paper, we will define it in general since the
definition is simple in the general case. Furthermore, we define the
curve by its ideal sheaf whose description we will use later on. We do
this in such a way that our nonreduced curves will be flat limits of
the reduced ones. Note that the restriction $d\geq e+2$ avoids trivial
cases where either the maps $X\ra A$ are constant or the cohomology
class of the image $\oX$ of $X$ is equal to the minimal class in which case
we know that $\oX$ does not deform out of the jacobian locus
\cite{matsusaka59}.

Let $W(L)\subset H^0 (L)$ be the vector subspace whose projectivization
is $L\subset |L|$. The underlying set of $X$ is the subset of $C^{ (e)
}$ where the elements of $W(L)$ are dependent. A scheme
structure can be defined on this set in the following way. Let
$D^e\subset C^{ (e) }\times C$ be the universal divisor and let $q_e$
and $p_e$ be the first and second projections from $C^{ (e) }\times C$
onto $C^{ (e) }$ and $C$ respectively. Then the global evaluation of
sections of $\cO_C (L)$ on divisors of degree $e$ is the map
\begin{equation}\label{detmap}
H^0 (L)\otimes\cO_{ C^{ (e) }}\lra V_L^e := {q_e}_* (p_e^*\cO_C (L) |_{ D^e })
\end{equation}
obtained by push-forward via $q_e$ from the evaluation map
\[
p_e^*\cO_C (L) \lra p_e^*\cO_C (L) |_{ D^e }.
\]
So $X$ is the locus where the evaluation map $W(L)\otimes\cO_{ C^{ (e)
}}\ra V_L^e$ has rank $\leq 1$. Therefore, since $X$ is of (the
expected) pure dimension $1$, by Eagon and Northcott
\cite{eagonnorthcott62} Theorem $2$ page $201$, there is an exact
sequence
\begin{equation}\label{eagon}
0\ra\Lambda^{e} {V^e_L}^*\otimes S^{ e-2 } W
(L)\ra\ldots\ra\Lambda^{4} {V^e_L}^*\otimes S^{ 2 } W
(L)\ra\Lambda^{3} {V^e_L}^*\otimes W (L)\ra\Lambda^{2}
{V^e_L}^*\ra\cI_{X / C^{ (e) } }\ra 0.
\end{equation}

Since our construction can be done globally in families, we see that
our non-reduced curves $X$ are indeed flat limits of reduced curves
$X$.

The natural morphism $C\ra\bP^1$ obtained from $L$
gives a morphism $X\ra\bP^1$ and, using the Hurwitz formula, one sees
immediately that $X$ has arithmetic genus
\[
g_X = -{d\choose e} +( g-1 +d ){d-2\choose e-1 } +1.
\]
This works at least when $X$ is smooth. When $X$ is not smooth, we
obtain the arithmetic genus by specialization from the smooth case.

\subsection{} Having defined $X$ in $C^{ (e) }$, we define the curve
$\oX$ that we are really interested in as its image in $A$ up to
translation. For this we first choose $g-e-1$ general points $p_1
,\ldots , p_{g-e-1}$ in $C$ and map $C^{ (e) }$ to $C^{ (g-1) }$ and
$A$ by the respective morphisms
\[
\begin{array}{rllllrll}
\phi_p : C^{ (e) } &\lra & C^{ (g-1) } & & & \psi_p : C^{ (e) } &\lra & A\\
D_e &\longmapsto & D_e +\sum_{ i=1 }^{ g-e-1 } p_i &
& & D_e &\longmapsto &\cO_C (D_e +\sum_{ i=1 }^{ g-e-1 } p_i).
\end{array}
\]
The first map is an embedding and the second map is a rational
resolution of its image which is a determinantal variety. The fibers
of $\psi_p$ are the complete linear systems in $C^{ (e) }$. Therefore,
in particular, if we let $W_e$ be the image of $C^{ (e) }$ in $A$,
then ${\psi_p}_* (\cO_{ C^{ (e) }} ) =\cO_{ W_e }$. We define $\oX$ to
be the curve whose ideal is ${\psi_p }_* (\cI_{ X/ C^{ (e) }})$. It
immediately follows that if $X$ is integral, then so is
$\oX$. Furthermore, it is the flat limit of general curves $\oX$ in
jacobians of general curves with $g^1_d$. Replacing $X$ by $X_{ d-e }
(L)$ if necessary, we will assume that $d\geq 2e$. We have the
following.

\begin{lemma}\label{h0e+11}
Suppose $e\leq g-2$ and $L$ contains reduced divisors. Then
\begin{enumerate}
\item there are divisors $D\in X_{ e+1 }(L)$ such that $h^0 (D) =1$,
\item assume $e\leq g-3$ and
\begin{enumerate}
\item either $d\geq 2e+2$,
\item or $d=2e+1$, $h^0 (L)\leq e$,
\item or $d=2e$, $h^0(L)\leq e-1$;
\end{enumerate}
then there are divisors $D\in X_{e+2} (L)$ such that $h^0 (D)=1$.
\end{enumerate}
\end{lemma}

\begin{proof}
If $d\geq 2g-2$, then a general divisor of $L$ is reduced and spans at
least a hyperplane in the canonical space of $C$. So we can choose a
subdivisor of degree $e+1$ (resp. $e+2$ if $e\leq g-3$) of it
which spans a linear subspace of dimension $e$ (resp. $e+1$) of
$|K|^*$ and hence
satisfies the lemma. If $d\leq 2g-3$, then, by Clifford's Theorem,
since $C$ is not hyperelliptic, we have $2 (h^0 (L) -1) < d$, hence $h^1
(L) < g-\frac{d}{2}\leq g-e$. So $h^1 (L)\leq g-e-1$ and a general
divisor of $L$ is reduced and spans a linear subspace of $|K|^*$ of
dimension at least $e$. Therefore it has a subdivisor of degree $e+1$
which spans a linear space of dimension $e$ and hence satisfies the
first part of the lemma. For the second part, the assumptions in
conjunction with Clifford's Theorem imply that $h^1 (L)\leq g-e-2$ and
an anlogous reasoning proves the second part.
\end{proof}

\section{The first order obstruction map}\label{sectobs}

\subsection{} From now on in the rest of the paper we shall
always assume that $X$ (hence $\oX$) is integral, i.e., reduced and
irreducible. It is immediate that the irreducibility of $X$ implies
that $L$ has no base points. Note that the converse to this is not
true as is easily seen by assuming that $C$ maps nontrivially to a
curve of positive genus and taking $L$ to be the inverse image of a
pencil on the curve of lower positive genus.

Recall that we also assume $d\geq 2e$, and, by Lemma
\ref{h0e+11}, a general $D_e\in X$ satisfies $h^0 (D_e ) =1$ so that
the map $X\ra\oX$ is birational.

\subsection{} Since $\oX$ is reduced, the obstructions to deformations
of $\oX$ with $A$ live in $Ext^1_{\oX} (\cI_{ \oX/A }/\cI_{ \oX/A }^2
,\cO_{\oX })$ (see
\cite{kollar96} Lemma 2.13 page 33 and Proposition 2.14 page 34).
We have the usual map
\begin{equation}\label{eqncotan}
\cI_{ \oX/A }/\cI_{ \oX/A }^2\lra\Omega^1_A |_{\oX }
\end{equation}
from which we obtain the map of exterior groups
\begin{equation}\label{eqnobs}
H^1 ( T_A |_{\oX} )\lra Ext^1 (\cI_{\oX/A }/\cI_{\oX/A }^2 ,\cO_{\oX} ).
\end{equation}
Composing this with restriction
\[
H^1 ( T_A )\lra H^1 ( T_A |_{\oX} ),
\]
we obtain the first order obstruction map
\[
\nu : H^1 (T_A )\ra Ext^1 (\cI_{\oX/A }/\cI_{\oX/A }^2 ,\cO_{\oX }).
\]
Given an infinitesimal deformation $\eta\in H^1 (T_A )$, the curve $X$
deforms with $A$ in the direction of $\eta$ if and only if $\nu (\eta
)= 0$.

\subsection{} The local to global spectral
sequence for the exterior sheaves of $\cI_{\oX/A }/\cI_{\oX/A }^2$ provides
the exact sequence
\[
0\lra H^1 ( N_{\oX/A })\lra Ext^1 (\cI_{\oX/A }/\cI_{\oX/A }^2 ,\cO_{\oX}
)\lra H^0 (\cE xt^1 (\cI_{\oX/A }/\cI_{\oX/A }^2 ,\cO_{\oX} )).
\]
The composition
\[
H^1 ( T_A )\lra Ext^1 (\cI_{\oX/A }/\cI_{\oX/A }^2 ,\cO_{\oX })\lra H^0 (\cE
xt^1 (\cI_{\oX/A }/\cI_{\oX/A }^2 ,\cO_{\oX} ))
\]
sends $\eta$ to the obstruction to deform $\oX$ with it locally. Since
$A$ is smooth, every deformation of $A$ is locally trivial and locally
$\oX$ deforms with it trivially. Therefore the image of $\eta$ by the
above composition is zero and the obstruction map $\nu$ factors
through $H^1 ( N_{\oX/A })$:
\[
\nu : H^1 (T_A )\lra H^1 ( N_{\oX/A }).
\]
Alternatively, the dual of the map (\ref{eqncotan}) gives us the map
of cohomology groups
\[
H^1 (T_A |_{\oX} )\lra H^1 ( N_{\oX/A })
\]
whose composition with the inclusion $H^1 ( N_{\oX/A })\lra Ext^1
(\cI_{\oX/A }/\cI_{\oX/A }^2 ,\cO_{\oX} )$ is (\ref{eqnobs}).

\subsection{} From the inclusion $\oX\subset W_e$, we obtain the
map
\[
N_{\oX/A }\lra N_{ W_e/A } |_{\oX}
\]
which gives us the map of cohomology groups
\[
H^1 (N_{\oX/A })\lra H^1 (N_{ W_e/A } |_{\oX}).
\]
We call $\nu_e$ the composition of this with $\nu$ and the pull-back
$H^1 ( N_{ W_e/A } |_{\oX} )\ra H^1 ( N_{ W_e/A } |_X)$ obtained from
the surjective morphism $X\ra\oX$:
\[
\nu_e : H^1 ( T_A )\lra H^1 ( N_{ W_e/A } |_X ).
\]
If $\nu_e (\eta )\neq 0$, then, a fortiori, $\nu (\eta )\neq 0$.

\subsection{} The choice of the polarization $\Theta$ provides an
isomorphism $H^1 (T_A )\cong H^1 (\cO_C )^{\otimes 2}$ via which the
algebraic (i.e. globally unobstructed) infinitesimal deformations with
which $\Theta$ deforms are identified with the elements of the subspace
$S^2 H^1 (\cO_C )\subset H^1 (\cO_C )^{\otimes 2}\cong H^1 (T_A )$. Via
this identification, the space of infinitesimal deformations of $(A
,\Theta)$ as a jacobian is naturally identified with $H^1 (T_C )\subset S^2
H^1 (\cO_C)$. The Serre dual of this last map is multiplication of sections
\[
S^2 H^0 ( K)\lra H^2 ( 2K )
\]
whose kernel is the space $I_2 (C)$ of degree $2$ forms vanishing on
the canonical image of $C$. To say that we consider an infinitesimal
deformation of $(A ,\Theta )$ out of the jacobian locus, means that we
consider $\eta\in S^2 H^1 (\cO_C )\setminus H^1 (T_C )$ which is
therefore equivalent to say that we consider $\eta\in S^2 H^1 (\cO_C
)$ such that there is $Q\in I_2 (C)$ with $(Q,\eta )\neq 0$. Here we
denote by
\[
(,) : S^2 H^0 (K)\otimes S^2 H^1 (\cO_C )\lra S^2 H^1 (K)\cong\bC
\]
the pairing obtained from Serre Duality.

\section{Translates of $\T$ containing $W_e$ and the obstruction
map}\label{sectrans}

To prove our main theorem we use translates of $\Theta$ which contain $W_e$.

\begin{lemma}
The subvariety $W_e$ is contained in a translate $\Theta_a$ of
$\Theta$ if and only if there exists $\sum q_i\in C^{ (g-e-1) }$ such
that $a =\sum p_i -\sum q_i$.
\end{lemma}
\begin{proof} For any points $q_1 ,\ldots ,q_{g-e-1}$ of $C$, the image of
$C^{ (e) }$ in $A$ by the corresponding map $\psi_q$ is contained in the
  divisor $\Theta_{\sum p_i -\sum q_i }$. Conversely, if $W_e$ is
  contained in a translate $\Theta_a$ of $\Theta$, then we have $h^0 ( D_e
  +\sum p_i -a) > 0$, for all $D_e\in C^{ (e) }$. Equivalently, for all
  $D_e\in C^{ (e) }$, we have $h^0 (K + a -\sum p_i -D_e ) > 0$, i.e., $h^0
  (K + a -\sum p_i )\geq e+1$ and $-a+\sum p_i$ is effective.
\end{proof}

\subsection{} Choose $a\in Pic^0 C$ such that $W_e\subset\Theta_a$ (i.e.,
  $a =\sum p_i -\sum q_i$ as above). Equivalently $W_e
-a\subset\Theta$. Let $\rho : C^{ (g-1 ) }\ra\Theta$ be the natural
morphism. Then (see \cite{green84} (1.20) p. 89) we have the exact
sequence
\begin{equation}\label{seqgr}
0\lra T_{ C^{ (g-1 )}}\lra\rho^* T_A\lra\cI_{ Z_{ g-1 }} (\Theta )\lra 0
\end{equation}
where, as in the Introduction, the scheme $Z_{ g-1 }$ is the locus
where the map $\rho$ fails to be an isomorphism. For the convenience
of the reader we mention that the scheme $Z_{ g-1 }$ is a
determinantal scheme of codimension $2$. If $g\geq 5$ or if $g=4$ and
$C$ has two distinct $g^1_3$'s, the scheme $Z_{ g-1}$ is reduced and
is the scheme-theoretical inverse image of the singular locus of
$\Theta$.

Combining sequence (\ref{seqgr}) with the tangent bundles sequences
for $C^{ (e) }_{\sum q_i}\subset C^{ (g-1) }$ and $W_e -a\subset A$,
we obtain the commutative diagram with exact rows and columns
\[
\begin{array}{ccccccccc}
 & & 0 & & 0 & & & & \\
 & & \downarrow & & \downarrow & & & & \\
 & & T_{ C^{ (e) } } & = & T_{ C^{ (e) } } & & & & \\
 & & \downarrow & &\downarrow & & & & \\
 & & T_{ C^{ (g-1) }} |_{ C^{ (e) }_{\sum q_i} } &\ra & T_A |_{ C^{ (e)
}_{\sum q_i} } &\ra & \cI_{ Z_{ g-1 }} (\Theta ) |_{ C^{ (e) }_{\sum q_i} }
&\ra & 0 \\
 & & \downarrow & &\downarrow & & \downarrow & & \\
 & & N_{  C^{ (e) }_{\sum q_i}/C^{ (g-1) }}  &\ra & N_{ {W_e -a
}/A } |_{ C^{ (e) }_{\sum q_i}} &\ra & \cO_{C^{ (e) }_{\sum q_i} }
(\Theta ) & \ra & 0 \\
 & & \downarrow & & & & & & \\
 & & 0 & & & & & & \\
\end{array}
\]
where the leftmost horizontal maps are injective if and only if $h^0
(\sum q_i ) =1$ and the map of sheaves
\[
T_A |_{ W_e -a}\lra N_{ {W_e -a }/A }
\]
fails to be surjective on the locus where $\psi_q$ fails to be
an embedding. In $C^{ (e) }$ this locus, which we call $Z_e$, is
locally defined, in the same way as $Z_{ g-1 }$, by the maximal minors
of the map $T_{ C^{ (e) }}\ra T_A |_{ C^{ (e) }}$. The support of
$Z_e$ is the subset of $C^{ (e) }$ parametrizing divisors $D_e$ with
$h^0 (D_e)\geq 2$. The right-hand bottom horizontal map is the
pull-back to $C^{ (e) }_{\sum q_i}$ of the map of normal sheaves
\[
N_{ {W_e -a }/A }\lra\cO_{ {W_e -a }} (\Theta) = N_{\Theta /A } |_{
W_e -a}
\]
whose image is the twist of a sheaf of ideals by $\cO_{ {W_e -a }}
(\Theta)$. As in the introduction, we let $Z_q$ be the subscheme of
$C^{ (e) }$ defined by the pull-back of this sheaf of ideals. Note
that because the sheaf of ideals contains $\cI_{ Z_{ g-1 }\cap C_{\sum
q_i}^{ (e) }}$, the subscheme $Z_q$ is contained in $Z_{ g-1 }\cap
C_{\sum q_i}^{ (e) }$. Furthermore, because $C^{ (e) }_{\sum q_i}\ra
W_e -a$ is an isomorphism outside $Z_e$, we have $Z_q\setminus Z_e =
Z_{ g-1 }\cap C_{\sum q_i}^{ (e) }\setminus Z_e$.

\subsection{Hypothesis} From now on we will assume that the set of $q$
for which $Z_q\cap X_{-a}\neq Z_{ g-1 }\cap X_{-a}$ has dimension at
most $g-e-3$ and we have chosen $q$ outside of this set, i.e., in such
a way that
\[
Z_q\cap X_{-a} = Z_{ g-1 }\cap X_{-a}.
\]
By Appendix \ref{generic} this will be the case for all $q\in C^{
(g-1-e) }$ on some open subset of the space of pairs $(C, L)$.

Therefore, restricting the previous diagram to $X_{ -a }$, we obtain
\[
\begin{array}{ccccccc}
T_{ C^{ (e) }_{\sum q_i} } |_{ X_{ -a }} & = & T_{ C^{ (e) }_{\sum q_i} } |_{
X_{ -a }} & & & & \\
\downarrow & &\downarrow & & & & \\
T_{ C^{ (g-1) }} |_{ X_{-a} } &\ra & T_A |_{ X_{-a} } &\ra &
\cI_{ Z_{ g-1 }} (\Theta ) |_{ X_{-a} } &\ra & 0 \\
\downarrow & &\downarrow & & || & & \\
N_{  C^{ (e) }_{\sum q_i}/C^{ (g-1) }} |_{ X_{ -a }} &\ra & N_{ {W_e -a
}/A } |_{ X_{ -a }} &\ra & \cI_{ Z_{ g-1 }} (\Theta ) |_{ X_{-a} }  &
\ra & 0 \\
\downarrow & & \downarrow & & & & \\
0 & & 0 & & & & \\
\end{array}
\]
So we have the commutative diagram
\[
\begin{array}{ccccc}\label{diagI}
S^2 H^1 (\cO_C )\subset H^1 ( T_A ) &\lra & H^1 (T_A |_{X_{-a}} )
&\lra & H^1 (\cI_{ Z_{ g-1 }} (\Theta ) |_{X_{-a}} ) \\
 || & & \downarrow & & || \\
S^2 H^1 (\cO_C )\subset H^1 ( T_A ) & \lra & H^1 ( N_{ {W_e -a }/A }
|_{ X_{ -a }} ) & \lra & H^1 (\cI_{ Z_{ g-1 } } (\Theta ) |_{X_{-a}} ).
\end{array}
\]
Translation by $a$ induces the identity on $H^1 (T_A )$ and isomorphisms
\[
\begin{split}
H^1 (T_A |_{X_{-a}} )\cong H^1 (T_A |_{X} )\\
H^1 ( N_{ {W_e -a }/A } |_{ X_{ -a }} )\cong H^1 ( N_{ W_e/A } |_X
)
\end{split}
\]
so that the kernel of
\[
\nu_e : S^2 H^1 (\cO_C )\lra H^1 ( N_{ W_e/A } |_X )
\]
is equal to the kernel of the map
\[
S^2 H^1 (\cO_C )\lra H^1 ( N_{ {W_e -a }/A } |_{ X_{ -a }} )
\]
obtained from $\nu_e$ by translation. Therefore the previous diagram
proves the following theorem.

\begin{theorem}\label{kernu}
The kernel of the map
\[
\nu_e : S^2 H^1 (\cO_C )\lra H^1 ( N_{ W_e/A } |_X )
\]
is contained in the kernel of the map obtained from the above
\[
S^2 H^1 (\cO_C )\lra H^1 (\cI_{ Z_{ g-1 }} (\Theta ) |_{X_{ -a }}
)
\]
for all $a$ such that $\Theta_{a}$ contains $W_e$ and $Z_q\cap X_{ -a
} = Z_{ g-1 }\cap X_{ -a }$.
\end{theorem}

A fortiori, the kernel of $\nu_e$ is contained in the kernel of the
composition
\[
S^2 H^1 (\cO_C )\lra H^1 (\cI_{ Z_{ g-1 }} (\Theta ) |_{X_{ -a }}
)\lra H^1 (\cI_{ Z_{ g-1 }\cap X_{ -a }} (\Theta )).
\]
We shall prove that for any $\eta\in S^2 H^1 (\cO_C
)\setminus H^1 ( T_C )$, there exists $a$ such that $\Theta_a$
contains $W_e$ and the image of $\eta$ in $H^1 (\cI_{ Z_{ g-1 }\cap
X_{-a}} (\Theta ))$ is nonzero unless
\begin{itemize}
\item either $e= h^0 ( L )$ and $d=2e$,
\item or $e = 2h^0 (L)$ and $d= 2e+1$.
\end{itemize}

\section{The kernel of the map $S^2 H^1 (\cO_C )\lra H^1 (\cI_{ Z_{
g-1 }\cap X_{ -a}} (\Theta ))$}\label{sectkertot}

\subsection{}\label{prepquad} The above map is equal to
the composition
\begin{equation}\label{comp2}
S^2 H^1 (\cO_C ) \lra H^1 (\cI_{ Z_{ g-1 } }(\Theta )) \lra H^1
(\cI_{ Z_{ g-1 } }(\Theta ) |_{ X_{ -a}}) \lra H^1 (\cI_{ Z_{ g-1
}\cap X_{ -a} } (\Theta )).
\end{equation}

From the usual exact sequence
\[
0\lra \cI_{ Z_{ g-1 } } (\Theta )\lra \cO_{ C^{ (e) }} (\Theta
)\lra\cO_{ Z_{ g-1 } } (\Theta )\lra 0
\]
we obtain the embedding
\[
H^0 (\cO_{ Z_{ g-1 }} (\Theta ))\inj H^1 (\cI_{ Z_{ g-1 }} (\Theta )).
\]
By \cite{green84} p. 95, the image of $S^2 H^1 (\cO_C )$ in $H^1
(\cI_{ Z_{ g-1 }} (\Theta ))$ is contained in $H^0 (\cO_{ Z_{ g-1 }}
(\Theta ))$. Using the commutative diagram with exact rows
\[
\begin{array}{ccccccccc}
0 &\lra &\cI_{ Z_{ g-1 }} (\Theta )&\lra &\cO_{ C^{ (g-1) }} (\Theta
)&\lra &\cO_{ Z_{ g-1 }} (\Theta )&\lra & 0\\
& &\downarrow & &\downarrow & &\downarrow & & \\
0 &\lra &\cI_{ Z_{ g-1 }\cap X_{-a}} (\Theta )&\lra &\cO_{X_{-a}} (\Theta
)&\lra &\cO_{ Z_{ g-1 }\cap X_{-a}} (\Theta )&\lra & 0,\\
\end{array}
\]
Composition (\ref{comp2}) is also equal to the composition
\[\begin{split}
S^2 H^1 (\cO_C )\lra H^0 (\cO_{ Z_{ g-1 }} (\Theta ))\lra 
H^0 (\cO_{ Z_{ g-1 }\cap {X_{-a}}} (\Theta ))\lra\\
\lra  H^1 (\cI_{ Z_{ g-1 }\cap {X_{-a}}} (\Theta )).\end{split}
\]
By \cite{green84} p. 95 the first map is the following
\[
\begin{array}{ccc}
S^2 H^1 (\cO_C ) &\lra & H^0 (\cO_{ Z_{ g-1 }} (\Theta ))\\
\sum a_{ij}\frac{\partial^2 }{\partial z_i\partial z_j} &\longmapsto
&\sum a_{ij}\frac{\partial^2\sigma }{\partial z_i\partial z_j} |_
{ Z_{ g-1 }},
\end{array}
\]
where $\{ z_i\}$ is a system of coordinates on $A$ and $\sigma$ is a
theta function with divisor of zeros equal to $\Theta$.
So we have the following description
\[
\begin{array}{ccccccc}
S^2 H^1 (\cO_C ) &\ra & H^0 (\cO_{ Z_{ g-1 }} (\Theta )) &\ra &
H^0 (\cO_{ Z_{ g-1 }\cap X_{-a}} (\Theta )) & \stackrel{coboundary}{\ra}
& H^1 (\cI_{ Z_{ g-1 }\cap X_{-a}} (\Theta )) \\
\sum a_{ij}\frac{\partial^2 }{\partial
z_i\partial z_j} &\mapsto &\sum a_{ij}\frac{\partial^2\sigma
}{\partial z_i\partial z_j} |_{ Z_{ g-1 }} &\mapsto &\sum
a_{ij}\frac{\partial^2\sigma }{\partial z_i\partial z_j} |_{ Z_{
g-1 }\cap X_{-a}} &\mapsto & ?
\end{array}
\]

\subsection{} We will investigate the kernel of the composition of the
first two maps $S^2 H^1 (\cO_C )\ra H^0 (\cO_{ Z_{ g-1 }\cap X_{-a}}
(\Theta ))$ and that of the coboundary map $H^0 (\cO_{ Z_{ g-1 }\cap
X_{-a}} (\Theta ))\ra H^1 (\cI_{ Z_{ g-1 }\cap X_{-a}} (\Theta ))$
separately. The kernel of $S^2 H^1 (\cO_C )\ra H^0 (\cO_{ Z_{ g-1
}\cap X_{-a}} (\Theta ))$ is contained in (with equality if and only
if $Z_{ g-1 }\cap X_{-a}$ is reduced) the annihilator of the quadrics
of rank $\leq 4$ which are the tangent cones to $\Theta$ at the points
of $Z_{ g-1 }\cap X_{-a}$.

\section{The kernel of the map $S^2 H^1 (\cO_C )\ra H^0 (\cO_{ Z_{ g-1
}\cap X_{-a}} (\Theta ))$}\label{sectker1}

\subsection{}\label{subsecker} An effective divisor $\sum_{ i=1 }^{
g-1-e } q_i\in C^{ (g-1-e)}$ gives the embedding of $X$ in $C^{ (g-1)
}$ defined by $D_e\mapsto D_e +\sum_{ i=1 }^{ g-1-e } q_i$. The union of
the images of these maps is a scheme, denoted $X + C^{ (g-1-e)
}\subset C^{ (g-1) }$ whose intersection with $Z_{ g-1}$ is the union
of the schemes $Z_{ g-1 }\cap X_{ -a }$ as $a=\sum p_i -\sum q_i$
varies.
To say that $\eta$ is in the kernel of the composition
\[
S^2 H^1 (\cO_C ) \lra H^0 (\cO_{ Z_{ g-1 }} (\Theta )) \lra
H^0 (\cO_{ Z_{ g-1 }\cap X_{-a}} (\Theta ))
\]
for every $a$, means that $\eta$ is in the kernel of the map
\[
S^2 H^1 (\cO_C ) \lra H^0 (\cO_{ Z_{ g-1 }} (\Theta )) \lra
H^0 (\cO_{ Z_{ g-1 }\cap X+ C^{ (g-1-e) }} (\Theta )).
\]
Since $X$ is reduced, $X+ C^{ (g-1-e) }$ is also the union of all
$C_E^{ (g-1-e) }$ with $E\in X$. So we see that the above is also
equivalent to $\eta$ being in the kernel of all the maps
\[
S^2 H^1 (\cO_C ) \lra H^0 (\cO_{ Z_{ g-1 }} (\Theta )) \lra
H^0 (\cO_{ Z_{ g-1 }\cap C_E^{ (g-1-e) }} (\Theta ))
\]
for all points $E\in X$.

\subsection{} We compute the kernel of the map
\[
S^2 H^1 (\cO_C ) \lra H^0 (\cO_{ Z_{ g-1 }} (\Theta )) \lra H^0
(\cO_{ Z_{ g-1 }\cap C_D^{ (f) }} (\Theta ))
\]
for $1\leq f\leq g-1$ and $D$ a divisor of degree $g-1-f$ such that
$h^0 (D) =1$. We shall later assume $D\in X$. Recall that by Lemma
\ref{h0e+11}, for a general $D\in X$, we have $h^0 (D) =1$.

Consider the exact sequence
\[
0\lra T_{ C^{ (g-1) }} |_{C_D^{ (f) }}\lra T_A |_{ C_D^{ (f) }}\lra
\cI_{Z_{ g-1 }} (\theta) |_{ C_D^{ (f) }}\lra 0
\]
from which it follows that the kernel of the map
\[
H^1 (T_A |_{ C^{ (f) }} ) = H^1 (T_A ) = H^1 (\cO_C )^{\otimes 2}\lra
H^1 (\cI_{Z_{ g-1 }} (\theta ) |_{ C_D^{ (f) }})
\]
is the image of $H^1 (T_{ C^{ (g-1) }} |_{ C_D^{ (f) }})$. We have
\[
\cI_{Z_{ g-1 }} (\theta ) |_{ C_D^{ (f) }}\cong\cI_{Z_{ g-1 }\cap
C_D^{ (f) }} (\theta )
\]
if and only if $Z_{ g-1 }\cap C_D^{ (f) }$ has codimension $2$ in
$C_D^{ (f) }$. In such a case the kernel of
\[
S^2 H^1 (\cO_C ) \lra H^0 (\cO_{ Z_{ g-1 }} (\Theta )) \lra H^0
(\cO_{ Z_{ g-1 }\cap C_D^{ (f) }} (\Theta ))
\]
is the intersection of the image of $H^1 (T_{ C^{ (g-1) }} |_{ C_D^{
(f) }})$ with $S^2 H^1 (\cO_C )$.

Denote by $V (D)$ the vector subspace of $H^1 (\cO_C ) = H^0 (\omega_C
)^*$ whose projectivization is $\langle D\rangle$.

\begin{theorem}
The image of $H^1 ( T_{ C^{ (g-1) }} |_{ C_D^{ (f) }}
)$ in $H^1 ( T_A |_{ C_D^{ (f) }} )$ is the span of $H^1 (T_C )\subset
S^2 H^1 (\cO_C )$ and $V (D)\otimes H^1 (\cO_C )$. The intersection
of this span with $S^2 H^1 (\cO_C )$ is the span of $H^1 (T_C )$ and
$S^2 V (D)$.
\end{theorem}

\begin{proof} We first determine the image of $H^1 (T_{ C^{ (g-1) }} |_{
C_D^{ (f) }})$ in $H^1 (T_A |_{ C_D^{ (f) }})$ and its intersection
with $S^2 H^1 (\cO_C )$. For this we use the diagram

\[
\begin{array}{ccccccc}
& 0 & & 0 & & & \\
& \downarrow & & \downarrow & & & \\
& T_{ C^{ (f) }} & =\hspace{-3pt}= & T_{ C^{ (f) }} & & & \\
& \downarrow & & \downarrow & & & \\
0\lra & T_{ C^{ (g-1) }} |_{C_D^{ (f) }} & \lra & T_A |_{ C_D^{ (f) }} &
\lra & \cI_{Z_{ g-1 }} (\theta) |_{ C_D^{ (f) }} & \lra 0 \\
& \downarrow & & \downarrow & & || & \\
0\lra & N_{ C_D^{ (f) }/ C^{ (g-1) }} & \lra & N_{ C_D^{ (f) }/A} & \lra &
\cI_{Z_{ g-1 }} (\theta ) |_{ C_D^{ (f) }} & \lra 0 \\
& \downarrow & & & & & \\
& 0 .& & & & &
\end{array}
\]
The map $H^1 ( C^{ (f) }, T_{C ^{ (f) }})\ra H^1 (C^{ (f) } , T_A |_{
C_D^{ (f) }})$ is injective, hence so is the map $H^1 ( C^{ (f) },
T_{C ^{ (f) }})\ra H^1 (C^{ (f) } , T_{C^{ (g-1) }} |_{ C_D^{ (f)
}})$. Therefore $H^0 ( C^{ (f) }, T_{ C^{ (g-1) }} |_{ C_D^{
(f) }} )\ra H^0 ( C^{ (f)} , N_{ C_D^{ (f) }/C^{ (g-1) }})$ is an
isomorphism.

Consider now the composition
\[
H^0 (T_{ C^{ (g-1)}} |_{ C_D^{ (f) }} )\otimes\cO_{ C^{ (f) }}\lra T_{
  C^{ (g-1) }} |_{ C_D^{ (f) }}\lra N_{ C_D^{ (f) }/ C^{ (g-1 )}}
\]
where the first map is evaluation. Then, because of the isomorphism
$H^0 ( C^{ (f) }, T_{ C^{ (g-1) }} |_{ C_D^{ (f) }} )\cong H^0 ( C^{
  (f)} , N_{ C_D^{ (f) }/C^{ (g-1) }})$, this composition can be
identified with evaluation of global sections
\[
H^0 ( C^{ (f)} , N_{ C_D^{ (f) }/C^{ (g-1) }})\otimes\cO_{ C^{ (f)
    }}\lra N_{ C_D^{ (f) }/C^{ (g-1) }}\; .
\]
From this we obtain the map
\[
H^0 ( C^{ (f)} , N_{ C_D^{ (f) }/C^{ (g-1) }})\otimes H^1 (\cO_{ C^{ (f)
    }} )\lra H^1 (N_{ C_D^{ (f) }/C^{ (g-1) }} ).
\]
Write $D =\sum_{ i=1 }^{ g-1-f } t_i$. Since $C_D^{ (f) }$ is the
complete intersection of the divisors $C_{ t_i }^{ (g-2) }$ in $C^{
(g-1) }$, its normal sheaf is isomorphic to $\oplus_{ i=1 }^f\cO_{ C^{
(f) }} ( C_{ t_i }^{ (f-1) })$. Therefore, using Appendix 6.1 in
\cite{I13}, the above map can be identified with
\[
\oplus_{ i=1 }^{ g-1-f } S^{ f-1 } H^0 (t_i )\otimes H^1 (\cO_C
)\lra\oplus_{i=1}^f S^{ f-1 } H^0 (t_i )\otimes H^1 (t_i)
\]
which is onto because each of the maps $H^1 (\cO_C )\ra H^1 (q_i )$ is
linear projection which is onto. Therefore, the composition
\[
H^0 (T_{ C^{ (g-1) }} |_{ C_D^{ (f) }} )\otimes H^1 (\cO_{ C^{ (f) }}
    )\lra H^1 (T_{ C^{ (g-1) }} |_{ C_D^{ (f) }} )\lra H^1 (N_{ C_D^{
    (f) }/C^{ (g-1) }} )
\]
is onto and hence so is
\[
H^1 (T_{ C^{ (g-1) }} |_{ C_D^{ (f) }} )\lra H^1 (N_{ C_D^{ (f) }/C^{
    (g-1) }} )\; .
\]
In conclusion we have the exact sequence
\[
0\lra H^1 ( T_{ C^{ (f) }} )\lra H^1 ( T_{ C^{ (g-1) }} |_{ C_D^{ (f)
}} )\lra H^1 ( N_{ C_D^{ (f) } / C^{ (g-1) }} )\lra 0
\]
and the image of $H^1 ( T_{ C^{ (g-1) }} |_{ C_D^{ (f) }} )$ in $H^1 (
T_A |_{ C_D^{ (f) }} )$ is the span of the images of $H^1 (T_C ) = H^1
( T_{ C^{ (f) }} )$ and $\oplus_{ i=1 }^{ g-1-f } S^{ f-1 } H^0 (t_i
)\otimes H^1 (\cO_C )$.

The image of $\oplus_{ i=1 }^{ g-1-f } S^f H^0 (t_i )$ in $H^0 (T_A )
= H^1 (\cO_C )$ is $V( D)$. Indeed, as we saw above, $H^0 ( T_{ C^{
(g-1) }} |_{ C_D^{ (f) }} )\cong H^0 (N_{ C_D^{ (f) }/C^{ (g-1) }} )
=\oplus_{ i=1 }^{ g-1-f } S^f H^0 (t_i )$ has dimension $g-1-f$. The
tangent space to $C^{ (g-1) }$ at $D_{ g-1 }\in C^{ (g-1) }$ can
canonically be identified with $\cO_{ D_{ g-1 }} (D_{ g-1 })$. For all
$D_{ g-1 }\in C_D^{ (f) }\subset C^{ (g-1) }$, we have $\cO_D
(D)\subset\cO_{ D_{ g-1 }} (D_{ g-1 })$. So $\cO_D (D)\subset H^0 (
T_{ C^{ (g-1) }} |_{ C_D^{ (f) }} )$ and the two spaces are equal
since they have the same dimension. The image of $\cO_D (D)$ in $H^0
(T_A ) = H^1 (\cO_C )$ by the differential of the map $\rho : C^{
(g-1) }\ra A$ is $V (D)$. So the image of $H^0 ( T_{ C^{ (g-1) }} |_{
C_D^{ (f) }} )$ in $H^1 (\cO_C )$ is $V (D)$. Therefore the image of
$H^0 ( T_{ C^{ (g-1) }} |_{ C_D^{ (f) }} )\otimes H^1 (\cO_C )
=\oplus_{ i=1 }^{ g-1-f } S^{ f-1 } H^0 (t_i )\otimes H^1 (\cO_C )$ in
$H^1 ( T_A |_{ C_D^{ (f) }} ) = H^1 (\cO_C )^{\otimes 2}$ is $V
(D)\otimes H^1 (\cO_C )$.
\end{proof}

Therefore, if $Z_{ g-1 }\cap C_D^{ (f) }$ has codimension $2$ in
$C_D^{ (f) }$, then the kernel of $S^2 H^1 (\cO_C )\ra H^0 (\cO_{ Z_{
g-1 }\cap C_D^{ (f) }})$ is the span of $H^1 (T_C )$ and $S^2 V (D)$.

Since we have assumed $h^0 (D) =1$, the codimension of $Z_{ g-1 }\cap
C_D^{ (f) }$ is at least $1$. If the codimension of $Z_{ g-1 }\cap
C_D^{ (f) }$ in $C^{ (f) }$ is $1$, then we have the exact sequence

\[
0 \lra \cK_Y \lra \cI_{ Z_{ g-1 }} (\T ) |_{ C^{ (f) }} \lra
\cI_{ Z_{ g-1 }\cap C^{ (f) }} (\T ) \lra 0
\]
where $Y$ is the maximal subscheme of $Z_{ g-1 }\cap C_D^{ (f) }$
supported on the union of its codimension $1$ components and $\cK_Y$
is the sheaf on $Y$ defined by the exact sequence. If
$f=1$, then $Y$ has dimension $0$, $h^1 (\cK_Y )=0$, $H^1 (\cI_{ Z_{
g-1 }} (\T ) |_{ C^{ (f) }} ) = H^1 (\cI_{ Z_{ g-1 }\cap C^{ (f) }}
(\T ))$ and the kernel of
\[
S^2 H^1 (\cO_C ) \lra H^0 (\cO_{ Z_{ g-1 }} (\Theta )) \lra H^0
(\cO_{ Z_{ g-1 }\cap C_D^{ (f) }} (\Theta ))\lra H^1 (\cI_{ Z_{ g-1
}\cap C^{ (f) }} (\T ))
\]
is again the intersection of the image of $H^1 (T_{ C^{ (g-1) }} |_{ C_D^{
(f) }})$ with $S^2 H^1 (\cO_C )$.

Suppose from now on that $2\leq f\leq g-1$. We have

\begin{lemma}\label{lemZcodim1}
The codimension of $Z_{ g-1 }\cap C_D^{ (f) }$ is $1$ only in the
following cases
\begin{enumerate}
\item The intersection $Z_{ g-1 }\cap C_D^{ (f) }$ contains
$C_t^{ (f-1) }$ for some $t\in C$. This happens if and only if
$\langle D\rangle\cap C$ contains $D+t$.
\item The restriction to $C$ of the projection from $\langle D\rangle$
is not birational to its image. Letting $C'$ be the normalization of
this image, the projection from $\langle D\rangle$ induces $\kappa :
C\ra C'$ of degree at least $2$. Given $\kappa$, there exist a finite
number of linear subspaces $L_i$ ($i= 1,\ldots ,l$) of $|\omega_C |^*$
such that any such $\langle D\rangle$ contains $L_i$ for some $i$.
Furthermore, $Z_{ g-1 }\cap C_D^{ (f) }$ contains the divisor $C^{
(f-2) } + X(\kappa )\subset C^{ (f) }$ where
\[
X (\kappa ) :=\{ D_2\in C^{ (2) } :\exists t\in C', h^0 (\kappa^* (t)
-D_2 ) > 0\}.
\]
\end{enumerate}
\end{lemma}

\begin{proof}
The first case is clear. Assume therefore that $Z_{ g-1
}\cap C_D^{ (f) }$ contains an irreducible divisor $\cF$ which is {\em
not} of the form $C_t^{ (f-1) }$. It is easily seen that this is
equivalent to the fact that for a general divisor $D_{ f-2 }\in C^{
(f-2) }$, the projection from $\langle D + D_{ f-2 }\rangle$ is not
birational on $C$. It first follows that the projection from $\langle
D\rangle$ is not birational on $C$. Indeed, if we call $C_1$ the image
of $C$ by the projection from $\langle D\rangle$, then, by the general
position theorem (\cite{ACGH} page 109), the projection of $C_1$ from
the span of a general effective divisor on it is always birational
unless the image of the projection is $\bP^1$ or a point. So, if the
projection of $C$ from $\langle D\rangle$ is birational, then so is
its projection from $\langle D + D_{ f-2 }\rangle$ for $D_{ f-2 }\in
C^{ (f-2) }$ general.

The general divisors $D_f\in\cF$ are of the form $D_{ f-2
} + D_2$ where $D_{ f-2 }\in C^{ (f-2) }$ is general and
$D_2\leq\kappa^* (t)$ for some $t\in C'$, i.e.,
\[
\cF = C^{ (f-2) } + X(\kappa ).
\]

To prove the assertion about the $L_i$, first suppose that the cover
$\kappa :C\ra C'$ is Galois and let $\{\sigma_1 ,\ldots ,\sigma_n\}$
be a set of generators for its Galois group. Then, since the
projection from $\langle D\rangle$ induces $\kappa$, the linear space
$\langle D\rangle$ is globally invariant under $\sigma_1 ,\ldots
,\sigma_n$ and $\sigma_1 ,\ldots ,\sigma_n$ induce the identity on
$|\omega_C |^* /\langle D\rangle$. Therefore, if we let $V$ be the
vector space whose projectivization is $|\omega_C |^* /\langle
D\rangle$, then, for all $i$, $\sigma_i$ has only one eigenvalue, say
$\lambda_i$ on $V$. Hence $\langle D\rangle$ contains the eigenspaces
of $\sigma_i$ for all its eigenvalues which are distinct from
$\lambda_i$. For each choice $\mu_1 ,\ldots ,\mu_n$ of eigenvalues of
$\sigma_1 ,\ldots ,\sigma_n$, we let $L (\mu_1 ,\ldots ,\mu_n)$ be the
smallest linear subspace of $|\omega_C |^*$ which, for all $i$,
contains all the eigenspaces of $\sigma_i$ distinct from $\mu_i$.
Assuming none of the $\sigma_i$ is the identity on $C$, all the $L(\mu_1
,\ldots ,\mu_n)$ are non-empty. So we see that $L(\lambda_1 ,\ldots
,\lambda_n )\subset\langle D\rangle$. It is immediate that a linear
subspace $L$ contains some $L (\mu_1 ,\ldots ,\mu_n)$ if and only if
the projection from $L$ factors through $\kappa : C\ra C'$. Therefore,
the $L (\mu_1 ,\ldots ,\mu_n)$ are the minimal subspaces $L$ of
$|\omega_C |^*$ such that the projection from $L$ factors through $\kappa :
C\ra C'$. This description shows that they only depend on $\kappa$ and
not the choice of the generating set $\{\sigma_1 ,\ldots
,\sigma_n\}$. We number them to obtain the subspaces $L_i$ ($i=
1,\ldots , l$) in the statement.

Now, if the cover $\kappa :C\ra C'$ is not Galois, it can be dominated
by a Galois cover: in other words, there exists a Galois cover
$\tkappa :\tC\ra C'$ which factors through $\kappa :C\ra C'$. The
induced map on the jacobians $J\tC\ra JC$ induces a projection
$|\omega_{\tC } |^*\ra |\omega_C |^*$ which, composed with the map
$|\omega_C |^*\ra |\omega_C |^* /\langle D\rangle$, induces $\tkappa$.
The subspaces $L_i$ are well-defined for $\tkappa$ and their images in
$|\omega_C |^*$ will give us the subspaces $L_i$ for $\kappa$.
\end{proof}

\begin{lemma}\label{lemDcodim2}
Suppose $e\leq g-3$. For $D\in X$ general, the intersection $Z_{ g-1
}\cap C_D^{ (g-1-e) }$ has codimension $2$ in $C_D^{ (g-1-e) }$.
\end{lemma}
\begin{proof}
We have to prove that neither of the two cases in Lemma
\ref{lemZcodim1} occur.

If, for $D$ general in $X$, $\langle D\rangle$ contains one of the
linear spaces $L_i$ from Lemma \ref{lemZcodim1}, then, since $X$ is
irreducible, for all $D\in X$, $\langle D\rangle\supset L_i$. Choose
now $D_{ e+1 }$ general in $X_{ e+1 } (L)$. Then, by Lemma \ref{h0e+11},
$h^0 (D_{ e+1 } )=1$ so that $\cap_{ D\leq D_{ e+1 }}\langle D\rangle
=\emptyset$ and it cannot contain any $L_i$.

Suppose now that, for $D$ general in $X$, $\langle D\rangle\cap
C\supset D+t_D$ for some $t_D\in C$. Then if we choose $D_{ e+1 }$ as
above and $D\leq D_{ e+1 }$ again, we see that $D+t_D\subset\langle
D_{ e+1 }\rangle$. If the $t_D$ are distinct for the distinct
subdivisors of $D_{ e+1 }$, then we obtain a contradiction by
Clifford's Theorem and the fact that $C$ is not hyperelliptic. Since
any two subdivisors of degree $e$ of $D_{ e+1 }$ have $e-1$ points in
common, we see then that there is a divisor $D_{ e-1 }\in X_{ e-1 }
(L)$ such that $\langle D_{ e-1 }\rangle\cap C\supset D_{ e-1 }+t$ for
some $t\in C$. Since $X_{ e-1 } (L)$ is also irreducible and our
choices of divisors were general, this is the case for all $D_{ e-1
}\in X_{ e-1 } (L)$. Repeating the argument with $e-1$ instead of $e$
and continuing, we arrive at a contradiction.
\end{proof}

\subsection{} Therefore, by what we saw above, for $D\in X$ general,
the kernel of $S^2 H^1 (\cO_C )\ra H^0 (\cO_{ Z_{ g-1 }\cap C_D^{ (f)
}})$ is the span of $H^1 (T_C )$ and $S^2 V (D)$. We have

\begin{lemma}\label{lemh0spanD}
Assume $e\leq g-3$ and
\begin{enumerate}
\item either $d\geq 2e+2$,
\item or $d=2e+1$, $h^0 (L)\leq e$,
\item or $d=2e$, $h^0(L)\leq e-1$.
\end{enumerate}
Then there is a reduced $D\in X_{ e+1 }(L)$ such that $h^0 (D) =1$ and
$\langle D\rangle\cap C =D$.
\end{lemma}

\begin{proof}
By Lemma \ref{h0e+11}, there is $E\in X_{ e+2 }(L)$ such that $h^0 (E)
=1$ and $E$ is reduced. We claim that there is $D\leq E$ with $\langle
D\rangle\cap C = D$. Suppose not so that the span of every subdivisor
of degree $e+1$ of $E$ contains an extra point of $C$. If two of these
points are equal, then we have a subdivisor of degree $e$ of $E$ whose
span contains an extra point of $C$ and this is not possible for all
such $E$ by the previous two Lemmas. So the span of $E$ contains
$e+2$ distinct extra points which is the number of subdivisors of
degree $e+1$ of $E$.  Therefore we have a divisor $E'$ of degree $e+2$
such that $\langle E+E'\rangle$ has dimension $e+1$. Hence $h^0 (E+E')
= 1+e+2 = e+3$. By Clifford's Theorem, since $C$ is not hyperelliptic,
this is possible only if $E+E'$ is a canonical divisor on $C$. In
particular, $e=g-3$.

Put $E = t_1 +\ldots t_{ g-1 }$ and $E' = s_1 +\ldots + s_{ g-1 }$,
the points being numbered in such a way that for all $j$, $s_j +\sum_{
i\neq j } t_i\leq\langle\sum_{ i\neq j } t_i\rangle\cap C$, i.e., $h^0
(s_j +\sum_{ i\neq j } t_i ) = 2$. Since $E+ E'$ is a canonical
divisor, we also have $h^0 (t_j +\sum_{ i\neq j } s_i ) =2$ for all
$j$, i.e., $t_j +\sum_{ i\neq j } s_i\leq\langle\sum_{ i\neq j }
s_i\rangle\cap C$. Choose a basis of $V(E) = V( E+ E')\subset H^1
(\cO_C )$ in which the coordinates of $t_j$ are $(0 ,\ldots
0,1,0\ldots ,0)$ where $1$ is in the $j$-th slot. Let $(a_{ j1 },\ldots
,a_{ j\; g-1 })$ be the coordinates of $s_j$. Then $a_{ jj } =0 $ for
all $j$. Take $j=1$. Then $a_{ 11 } =0$ and the condition
$t_1\in\langle\sum_{ i=2 }^{ g-1 } s_i\rangle$ means there are scalars
$\lambda_2 ,\ldots ,\lambda_{ g-1 }$ such that
\[
\begin{array}{l}
1 =\sum_{ i=2 }^{ g-1 }\lambda_i a_{ i1 }\\
0 =\sum_{ i=2 }^{ g-1 }\lambda_i a_{ ik }\hbox{ for all }k\geq 2
\end{array}
\]
Since $E+ E'$ is a canonical divisor and $h^0 (E) =1$, we also have
$h^0 (E') =1$. Therefore the $s_j$ are linearly independant and, a
fortiori, the minor $| a_{ jk } |_{\substack{2\leq j\leq g-1 \\ 2\leq
k\leq g-1}}$ is not zero and the condition $0 =\sum_{ i=2 }^{ g-1
}\lambda_i a_{ ik }\hbox{ for all }k\geq 2$ implies $\lambda_i =0$ for
all $i$. Then the condition $1 =\sum_{ i=2 }^{ g-1 }\lambda_i a_{ i1
}$ gives a contradiction.
\end{proof}

\begin{lemma}
Suppose $1\leq e\leq g-3$ and $D := t_1 +\ldots + t_{ e+1 }$ is a
reduced divisor such that $h^0 (D) =1$ and $\langle D\rangle\cap C =
D$. Put $E_i := D - t_i$. Then
\[
\bigcap_{ i=1 }^{ e+1 } \langle H^1 (T_C ), S^2 V (E_i)\rangle = H^1 (T_C ).
\]
\end{lemma}

\begin{proof}
We proceed by induction on $e$. For $e=1$, we have $E_i = t_{ 2-i}$,
\[
S^2 V(E_i )\subset H^1 (T_C),
\]
and the result is trivially true. Suppose $e\geq 2$ and the
result holds for $e-1$. Let us rewrite
\[
\bigcap_{ i=1 }^{ e+1 } \langle H^1 (T_C ), S^2 V (E_i)\rangle =
\bigcap_{ i=1 }^{ e }\left( \langle H^1 (T_C ), S^2 V (E_i)\rangle
\cap\langle H^1 (T_C ), S^2 V (E_{ e+1 })\rangle\right).
\]
We will prove that
\[
\langle H^1 (T_C ), S^2 V (E_i)\rangle\cap\langle H^1 (T_C ), S^2 V
(E_{ e+1 })\rangle =\langle H^1 (T_C ), S^2 V (E'_i)\rangle
\]
where $E'_i := D- t_i - t_{ e+1 } = E_{ e+1 } - t_i$. Then replacing
$D$ with $E_{ e+1 }$ we are reduced to the statement for $e-1$.

Dually, we will prove that the annihilators in $S^2 H^0 (\omega_C )$
of the two spaces are equal. The annihilator of $H^1 (T_C)$ is $I_2
(C)$. That of $\langle H^1 (T_C) , S^2 V(E)\rangle$ for any divisor
$E$ is the space $I_2 (C,E)$ of homogeneous degree $2$ forms vanishing
on $C$ and the linear span $\langle E\rangle$ of $E$ in $|\omega_C
|^*$. The statement we need to prove has now become
\[
I_2 (C, D- t_i -t_{ e+1 }) = I_2 (C , D- t_i )+ I_2 (C, D - t_{ e+1 }).
\]
Choose $g-3-e$ general points $t_{ e+2 },\ldots t_{ g-2 }$ on $C$.
Then
\[
\langle\sum_{ i=1 }^{ g-2 } t_i\rangle\cap C =\sum_{ i=1 }^{ g-2 } t_i
\]
and the $t_i $ are linearly independent and distinct. In particular,
the $t_i$ impose independent conditions on quadrics.

We claim that the restriction map
\[
I_2 (C)\lra S^2 V(\sum_{ i=1 }^{ g-2 } t_i)^*
\]
induces an isomorphism between $I_2 (C)$ and the homogeneous degree
$2$ forms on $V(\sum_{ i=1 }^{ g-2 } t_i)$ vanishing at the points
$t_i$. These two spaces have the same dimension so it is sufficient to
prove that the restriction map is injective, i.e., no quadric in
$|\omega_C |^*$ containing the canonical curve $C$ contains
$\langle\sum_{ i=1 }^{ g-2 } t_i\rangle$. Since $\langle\sum_{ i=1 }^{
g-2 } t_i\rangle$ has codimension $2$, if a quadric contains it, then
the quadric has rank $\leq 4$. Then $\langle\sum_{ i=1 }^{ g-2 }
t_i\rangle$ is a member of a ruling of the quadric and by
\cite{andreottimayer67} cuts a divisor of a $g^1_d$ on $C$. This,
however, is not possible by our assumptions on the $t_i$.

It first follows from this that
\[\begin{split}
\hbox{dim} I_2 (C, D- t_{ e+1 }) = \hbox{dim} I_2 (C, D- t_{ i }) =
\hbox{dim} I_2 (C) - {e\choose 2} = {g-2\choose 2} - {e\choose 2}\\
\hbox{dim} I_2 (C, D- t_{ e+1 } - t_i) =\hbox{dim} I_2 (C) -
{e-1\choose 2} = {g-2\choose 2} - {e-1\choose 2}.
\end{split}\]
So to prove our claim we need to prove that
\[\begin{split}
\hbox{dim}\left( I_2 (C, D- t_{ e+1 })\cap I_2 (C, D- t_{ i })\right) =
2\left({g-2\choose 2} - {e\choose 2}\right) -\left({g-2\choose 2} -
{e-1\choose 2}\right) \\
= {g-2\choose 2} - {e\choose 2} - (e-1)
\end{split}\]
This is easily seen to be true from our assumptions on the $t_i$.
\end{proof}

\subsection{} So far it follows from our results above that the intersection
of the kernels of the maps
\[
S^2 H^1 (\cO_C ) \lra H^0 (\cO_{ Z_{ g-1 }} (\Theta )) \lra
H^0 (\cO_{ Z_{ g-1 }\cap C_E^{ (g-1-e) }} (\Theta ))
\]
as $E$ varies in $X$ is $H^1 (T_C)$. Therefore (see \ref{subsecker})
for a given $\eta\not\in H^1 (T_C)$, there exists $a = \sum p_i -\sum
q_i$ such that $\eta$ is {\em not} in the kernel of the map
\[
S^2 H^1 (\cO_C ) \lra H^0 (\cO_{ Z_{ g-1 }} (\Theta )) \lra
H^0 (\cO_{ Z_{ g-1 }\cap X_{-a}} (\Theta )).
\]

\section{The kernel of the map $H^0 (\cO_{ Z_{ g-1 }\cap X_{-a}}
(\Theta ))\ra H^1 (\cI_{ Z_{ g-1 }\cap X_{ -a }} (\Theta ))$}\label{sectker2}

We continue with the analysis of the kernel of the coboundary map
\[
H^0 (\cO_{ Z_{ g-1 }\cap X_{-a}} (\Theta ))\lra H^1 (\cI_{ Z_{ g-1
}\cap X_{ -a }} (\Theta ))
\]
and see that in fact we can choose our $a$ above also in such a way that
this map is injective.

\subsection{}\label{sectdefZ} Let $\tZ (X)\subset C^{ (g-1-e) }\times
X$ be the closure of the subvariety parametrizing pairs $(\sum q_i ,
D_e )$ such that $h^0 (\sum q_i ) = h^0 ( D_e) = 1$ and $h^0 (\sum q_i
+ D_e )\geq 2$. Let $Z(X)\subset C^{ (g-1-e) }$ be the image of $\tZ
(X)$ by the first projection.
\begin{lemma}\label{lemZnotempty}
The varieties $\tZ (X)$ and $Z(X)$ are not empty.
\end{lemma}
\begin{proof}
Choose a general $D_e\in X$ so that we have $h^0 (D_e) =1$ (see Lemma
\ref{h0e+11} or \ref{lemh0spanD}). If, for
all $\sum q_i\in C^{ (g-1-e) }$ with $h^0 (\sum q_i + D_e )\geq 2$, we
have $h^0 (\sum q_i )\geq 2$, then, for some $r\geq 1$, the dimension
of $W^r_{ g-1-e }$ is at least $g-1-e-r-2 = g-e-r-3$. By
\cite{mumford74} pp. 348-350, this can only be the case if $r=1$ and
either $C$ is trigonal, bielliptic or a smooth plane quintic.

In the trigonal case $W^1_{ g-1-e } = g^1_3 + C^{ (g-e-4) }$. For a
point $t$ of $D_e$, the divisor $\sum q_i = g^1_3 -t + D_{ g-3-e }$ with
$D_{ g-3-e }\in C^{ (g-1-3) }$ general satisfies $h^0 (\sum q_i ) = 1$
and $h^0 (\sum q_i + D_e )\geq 2$.

In the bielliptic case, if $\pi : C\ra E$ is the bielliptic cover,
then $W^1_{ g-1-e } =\pi^* W^1_2 (E) + C^{ (g-e-5) }$. For two
distinct points $s$ and $t$ of $D_e$, the divisor $\iota s +\iota t +
D_{ g-3-e }$ with $D_{ g-3-e }\in C^{ (g-1-3) }$ general and $\iota$
the bielliptic involution satisfies $h^0 (\sum q_i ) = 1$ and $h^0
(\sum q_i + D_e )\geq 2$.

In the case of the smooth plane quintic, we have $g-1-e = 4 = g-2$. So
$e=1$ which is excluded.
\end{proof}

\begin{lemma}\label{lemcobnotinj}
Suppose that for $\sum q_i\in Z(X)$ with $h^0 (\sum q_i) =1$ the
coboundary map
\[
H^0 (\cO_{ Z_{ g-1 }\cap X} (\Theta_a )) \lra H^1 (\cI_{ Z_{ g-1
}\cap X} (\Theta_a ))
\]
is not injective, then
\[
H^0 (K -\sum q_i -L )\neq 0.
\]
\end{lemma}
\begin{proof}
Using the exact sequence
\[
0\lra H^0 (\cI_{ Z_{ g-1 }\cap X} (\Theta_a )) \lra H^0 (\cO_X
(\Theta_a )) \lra H^0 (\cO_{ Z_{ g-1 }\cap X} (\Theta_a )) \lra H^1
(\cI_{ Z_{ g-1 }\cap X} (\Theta_a )),
\]
we need to understand the sections of $\cO_X (\Theta_a )$ which vanish
on $Z_{ g-1 }\cap X$. For this, we use the embedding of $X$ in
$C^{(e)}$:
\[
0\lra\cI_{ X/C^{ (e) }} (\Theta_a )\lra\cO_{ C^{ (e) }} (\Theta_a )\lra\cO_X
(\Theta_a )\lra 0.
\]
By Appendix 6.1 in \cite{I13} this gives the exact sequence of cohomology
\[\begin{split}
0\lra H^0 (\cI_{ X/C^{ (e) }} (\Theta_a ))\lra\wedge^e H^0 (C, K
-\sum_{ i=1 }^{ g-e-1 } q_i)\lra H^0 (X,\Theta_a )
\lra H^1 (\cI_{ X/C^{ (e) }} (\Theta_a ) ).
\end{split}
\]
By Appendix 6.2 in \cite{I13} the elements of $H^0 ( C^{ (e)
},\Theta_a ) =\wedge^e H^0 (C, K -\sum_{ i=1 }^{ g-e-1 } q_i)$ all
vanish on $Z_{ g-1 }\cap C^{ (e) }$, hence they also vanish on $Z_{
g-1}\cap X$. So if the coboundary map is not injective, then there
must be elements of $H^0 (X,\Theta_a )$ which are not restrictions of
elements of $H^0 ( C^{ (e) },\Theta_a )$.  In particular, we must have
$H^1 (\cI_{ X/C^{ (e) }} (\Theta_a ) )\neq 0$. By sequence
(\ref{eagon}), this implies that there is an integer $j\in\{ 1,\ldots
,e\}$ such that $H^{ j-1 } (\Lambda^j V_L^{ e*} (\Theta_a ))\neq
0$. Equivalently, $H^{ e-j+1 } (\omega_{ C^{ (e) }} (-\Theta_a
)\otimes\Lambda^j V_L^e )\neq 0$. By Appendix \ref{cohVL}, since
(\cite{I13} Appendix)
\[
\pi_e^*\cO_{ C^{ (e) }} (\Theta_a )\cong
p_1^*\cO_C ( K-q )\otimes\ldots\otimes p_e^*\cO_C ( K-q ) (-\sum_{
1\leq k<l\leq e }\Delta_{ k,l }),
\]
this implies that
\[
H^{ e-j+1 } ( p_1^*\cO_C (L +q )\otimes\ldots\otimes p_j^*\cO_C (L +q)\otimes
p_{ j+1 }^*\cO_C (q)\otimes\ldots p_e^*\cO_C (q)
(-\sum_{ 1\leq k< l\leq j }\Delta_{ k,l }))^{\Sg_j\times\Sg_{ e-j }}\neq 0
\]
As in the Appendix of \cite{I13} the above cohomology group is equal
to the group of elements of
\[
H^{ e-j+1 } ( p_1^*\cO_C (L +q )\otimes\ldots\otimes p_j^*\cO_C (L +q)\otimes
p_{ j+1 }^*\cO_C (q)\otimes\ldots p_e^*\cO_C (q))
\]
anti-invariant under the action of $\Sg_j$ and invariant under the
action of $\Sg_{ e-j }$. Therefore its non-vanishing implies the
non-vanishing of $H^1 (L +q) = H^1 (L +\sum q_i ) = H^0 ( K-\sum q_i
-L)^*$.
\end{proof}
\subsection{} By Lemma \ref{lemZnotempty}, the variety $\tZ (X)$ is not
empty. Therefore $\tZ(X)$ has dimension at least $g-e-2$. By Lemma
\ref{lemtZZfin} below this implies that $Z(X)$ also has dimension $\geq
g-e-2$. Therefore the hypothesis and hence the conclusion of Lemma
\ref{lemcobnotinj} hold for a $(g-e-2)$-dimensional family of $\sum_{ i=1 }^{
g-e-1 } q_i$. This implies $h^0 (K-L )\geq g-e-1$. Therefore, by
Clifford's Lemma, since $C$ is not hyperelliptic, we have $2( g-e-1-1)
< 2g-2-d$ or $d\leq 2e+1$. If $d=2e+1$, then $h^0 (L) =e+1$ and $C$
has Clifford index $1$ if $e\leq g-3$.

\begin{lemma}\label{lemtZZfin}
The projection $\tZ (X)\ra Z(X)$ is generically finite.
\end{lemma}

\begin{proof}
If not, then some component of $\tZ (X)$ maps with one-dimensional
fibers into $Z (X)$. Since $X$ is integral, these one-dimensional
fibers are all isomorphic to $X$. Hence there is a
$(g-e-3)$-dimensional family of divisors $q=\sum q_i$ such that $h^0
(q )= 1$ and, for every $D_e\in X$, $h^0 (q+ D_e )\geq 2$. Let $e'$ be
the largest integer such that for a general such $q$ and a general
$D_{ e' }\in X_{ e' } (L)$, we have $h^0 (q + D_{ e'} ) =1$. It is
immediate that $1\leq e'\leq e-1$. For a fixed general such $D_{ e'
}$, let $D$ be the divisor of $L$ containing $D_{e'}$ and let $t$ be a
point of $D - D_{ e'}$. Then $h^0 ( q + D_{ e' } +t )=
2$. Equivalently $h^0 (K -q -D_{ e' } - t) = h^0 (K -q -D_{ e'})$,
i.e., $t$ is a base point of the linear system $|K - q - D_{ e'
}|$. Since $D_{ e' }$ is general, so is $D$, hence $D$ is reduced and,
furthermore, it has no points in common with $q$. It follows that all
of $D -D_{ e'}$ is contained in the base locus of $|K - q - D_{ e'
}|$. By Riemann-Roch and Serre duality we see that this implies $h^0
(q+D )= 1+d-e'$. So $C$ has a family of dimension $g-e-3$ of linear
systems of degree $g-1-e+d$ and dimension $d-e'$. Since $C$ is not
hyperelliptic, it follows from \cite{mumford74} pp. 348-350 that
$g-e-3\leq g-1-e+d -2(d-e') -1$, i.e., $d\leq 2e' +1\leq 2e-1$ which
contradicts the hypothesis $d\geq 2e$.
\end{proof}

\section{The consequences of Theorem \ref{mainthm}}\label{sectconsq}

If $C$ has a $g^{ e-1 }_{ 2e }$ (resp. $g^e_{ 2e+1 }$), then $C$ has
Clifford index $2$ (resp. $1$). By a result of Martens
(\cite{martens80} Satz 4 page 80), if $C$ is non-bielliptic, has
Clifford index $2$ (resp. $1$) and genus at least $10$ (resp. $8$),
then $C$ has no $g^{e-1}_{ 2e}$ for $4\leq e\leq g-5$ (resp. no $g^e_{
2e+1 }$ for $3\leq e\leq g-5$). The cases of low genus are easily
analyzed \cite{I16} and we see that the cases in which $X$ might
deform with $JC$ out of the jacobian locus are
\begin{enumerate}
\item $C$ any curve, $e= g-2$,
\item $C$ with a $g^1_4$ and hence also a $g^{g-4}_{ 2g-6} = |K_C -
g^1_4|$, $e=2$ or $g-3$,
\item $C$ with a $g^2_6$, $e=3$ or $g-4$,
\item $C$ with a $g^1_3$, $e= g-4$ or $g-3$,
\item $C$ with a $g^2_5$, $e=2$ or $g-4$,
\item $C$ bielliptic, $2\leq e\leq g-2$.
\end{enumerate}

\section{Appendix}

\subsection{The cohomology of $\cI_{ X/C^{ (e) }}\otimes N$}\label{cohVL}

Here we introduce a method for computing the cohomology of $\cI_{
X/C^{ (e) }}\otimes N$ where $N$ is a locally free sheaf on $C^{ (e)
}$. One way to approach this calculation is to compute the
cohomologies of the pieces $\Lambda^j V_L^{ e* }\otimes S^{ j-2 }
W(L)\otimes N$ of the resolution (\ref{eagon}) of $\cI_{ X/C^{ (e)
}}\otimes N$. Or equivalently, the cohomologies of the sheaves
$\omega_{ C^{ (e) }}\otimes N^*\otimes\Lambda^j V_L^e$. Recall that
\[
V_L^e = {q_e}_* (p_e^*\cO_C (L) |_{ D^e })
\]
where $D^e\subset C^{ (e) }\times C$ is the universal divisor and
$q_e, p_e$ are the first and second projections of $C^{ (e) }\times C$
onto its two factors. On this model, for $1\leq j\leq e$, let
\[
Y^{ e,j }\subset C^{ (e) }\times C^{ (j) }
\]
be the universal subvariety, i.e.,
\[
Y^{ e,j } :=\left\{ (D_e , D_j )\in C^{ (e) }\times C^{ (j) } : D_e\geq
D_j\right\},
\]
and let $q_{ e,j }$ and $p_{ e,j }$ be the first and second
projections of $C^{ (e) }\times C^{ (j) }$ onto its two factors. Then
a moment of reflexion will convince the reader that
\[
\Lambda^j V_L^e = {q_{ e,j }}_*\left(\left( p_{ e,j }^*\cL'_{ L,j
}\right) |_{ Y^{ e,j }}\right)
\]
where, as in \cite{I13}, $\cL'_{ L,j }$ is the sheaf on $C^{ (j) }$
whose inverse image on $C^j$ is $p_1^*\cO_C (L)\otimes\ldots
p_j^*\cO_C (L)\left( -\sum_{ 1\leq k<l\leq j }\Delta_{ k,l }\right)$. So
\[
H^k (\omega_{ C^{ (e) }}\otimes N^*\otimes\Lambda^j V_L^e ) = H^k
(\omega_{ C^{ (e) }}\otimes N^*\otimes {q_{ e,j }}_*\left(\left( p_{
e,j }^*\cL'_{ L,j }\right) |_{ Y^{ e,j }}\right) )
\]
is a graded piece of the filtration of
\[
H^k (q_{ e,j }^*\left(\omega_{ C^{ (e) }}\otimes N^*\right)\otimes p_{
e,j }^*\cL'_{ L,j } |_{ Y^{ e,j }})
\]
induced by the Leray Spectral sequence of the fibration $q_{ e,j } :
C^{ (e) }\times C^{ (j) }\ra C^{ (e) }$. The morphism
\[\begin{array}{ccc}
C^e & \surj & Y^{ e,j }\\
(s_1 ,\ldots , s_e) & \longmapsto & (s_1 +\ldots + s_e , s_1 +\ldots +
s_j)
\end{array}\]
shows that $Y^{ e,j}$ is the quotient of $C^e$ by the action of
$\Sg_j\times\Sg_{ e-j }$ which permutes the first $j$ points and the
last $e-j$ points. Therefore
\[\begin{split}
H^k (q_{ e,j }^*\left(\omega_{ C^{ (e) }}\otimes N^*\right)\otimes p_{
e,j }^*\cL'_{ L,j } |_{ Y^{ e,j }}) = H^k (\pi_{ e,j }^*\left(q_{ e,j
}^*\left(\omega_{ C^{ (e) }}\otimes N^*\right)\otimes p_{ e,j
}^*\cL'_{ L,j } |_{ Y^{ e,j }}\right))^{\Sg_j\times\Sg_{ e-j }} \\
= H^k ( p_1^*\cO_C (K +L)\otimes\ldots\otimes p_j^*\cO_C (K +L)\otimes
p_{ j+1 }^*\cO_C (K)\otimes\ldots\otimes p_e^*\cO_C (K)\otimes\pi_e^* N^*\\
(-\sum_{ 1\leq k< l\leq j }\Delta_{ k,l } -\sum_{ 1\leq k< l\leq e
}\Delta_{ k,l }))^{\Sg_j\times\Sg_{ e-j }}.
\end{split}\]

\subsection{}\label{generic}

The purpose of this section is to show that when $e\leq g-3$, for a
sufficiently general pair $(C, L)$ the curve $X$ obtained satisfies
the condition
\[
X\cap Z_{ g-1 } = X\cap Z_q
\]
for all $\sum q_i\in C^{ (g-1-e) }$. Since the schemes $Z_{ g-1 }\cap
(W_e -a)$ and $Z_q$ are equal outside $Z_e$, this will follow if we show
that $X\cap Z_e =\emptyset$.

So we shall prove that for $(C, L)$ sufficiently general, there is
no subdivisor $D$ of degree $e$ of a divisor of $L$ such that $h^0
(D)\geq 2$. If $d\leq g+1$, then for a sufficietly general $(C, L)$,
$h^0 (L) =2$ and since we can suppose $L$ base-point-free, the
assertion is true.

Suppose therefore that $d\geq g+2$. In such a case, supposing $(C, L)$
general also means that $C$ is general, since a general curve has
$g^1_d$'s. If $e <\frac{g+2}{2}$, then by Brill-Noether theory a
general curve does not have any $g^1_e$ and the assertion is true.

So we suppose $e\geq\frac{g+2}{2}$ (which then implies $d\geq g+2$). In
this case we prove that in a general linear system $G$ of degree
$d\geq g+2$ on the general curve $C$, the family $M_G$ of divisors of
the form $D_e + D_{ d-e }$ with $h^0 (D_e )\geq 2$ has codimension at
least $2$. Then a general pencil $L$ in $G$ will not intersect $M_G$.
It suffices to show that the union $M :=\cup_{ deg(G) =d} M_G$ has
dimension at most $d-2$ since $\cup_{ deg(G) =d} G = C^{ (d) }$ has
dimension $d$. We can rewrite $M = C^{ (e) }_1 + C^{ (d-e) } :=\cup_{
D\in C^{ (e) }_1} C^{ (d-e) }_D$. Since $C$ is general, by
Brill=Noether theory
\[
\hbox{dim} W^1_e = g-2( g-e+1 ) = 2e -g -2.
\]
Hence dim$C^{ (e) }_1 = 2e -g -1$ and dim$M = 2e -g -1 + d-e = d+
e-g-1$. We have $e-g-1\leq -2$ which concludes our proof.


\providecommand{\bysame}{\leavevmode\hbox to3em{\hrulefill}\thinspace}

\end{document}